\documentclass[10pt,leqno]{amsart}

\usepackage{soul}
\usepackage {yfonts}
\usepackage [T1] {fontenc}

\usepackage{xcolor}
\usepackage{textcomp}

\usepackage{enumitem}

\numberwithin{equation}{section}
\newtheorem{theorem}{Theorem}[section]
\newtheorem{prop}[theorem]{Proposition}
\newtheorem{lem}[theorem]{Lemma}

\newtheorem{rem}[theorem]{Remark}

\newcommand{\R}{\mathbb{R}}

\newcommand{\N}{\mathbb{N}}

\newcommand{\F}{\mathcal{F}}
\theoremstyle{plain}

\def\pn{\hfil\par\noindent}
\def\be{\begin{enumerate}}
\def\ds{\displaystyle }
\def\ee{\end{enumerate}}

\def\beqq{\begin{eqnarray*}}
\def\eeqq{\end{eqnarray*}}

\def\buildo#1\over#2{\mathrel{\mathop{\null#2}\limits^{#1}}}
\def\buildu#1\under#2{\mathrel{\mathop{\null#2}\limits_{#1}}}

\title{Long Time Decay of Leray Solution of 3D-NSE With Exponential Damping}
\author{Mongi Blel}
   \address{King Saud University, College of Sciences, Department of Mathematics,   Kingdom of Saudi Arabia}
   \email{mblel@ksu.edu.sa, jamelbenameur@gmail.com}
\author{Jamel Benameur}

\date{\today}

\subjclass[MSC 2020]{Primary  35-XX, 35Q30, 76D05, 76N10}
\keywords{Navier-Stokes Equations, Friedrich method, global weak solution}

\begin{document}

\maketitle

\begin{abstract}
We study the uniqueness, the continuity in $L^2$ and the large
time decay for the Leray solutions of the $3D$ incompressible
Navier-Stokes equations with nonlinear exponential
damping term $a (e^{b |u|^{\bf 4}}-1)u$, ($a,b>0$).
\end{abstract}

\section{\bf Introduction}\ \\

 In this paper, we investigate the questions of the  existence, uniqueness and asymptotic study  of global weak solution to the modified incompressible Navier-Stokes equations in $\R^3$

\begin{equation}\label{$S_2$}
 \left\{ \begin{matrix}
     \partial_t u
 -\nu\Delta u+ u.\nabla u  +a (e^{b |u|^4}-1)u = -\nabla p \hfill&\hbox{ in } \mathbb R^+\times \mathbb R^3\\
     {\rm div}\, u = 0 \hfill&\hbox{ in } \mathbb R^+\times \mathbb R^3\\
    u(0,x) =u^0(x) \;\;\hfill&\hbox{ in }\mathbb R^3,\\
    a,b>0\hfill&
\end{matrix}\right. \tag{$S$}
\end{equation}
 where $u=u(t,x)=(u_1,u_2,u_3)$, $p=p(t,x)$ denote respectively the unknown velocity and the unknown pressure of the fluid at the point $(t,x)\in \mathbb R^+\times \mathbb R^3$,  the viscosity of fluid $\nu>0$ and $u^0=(u_1^0(x),u_2^0(x),u_3^0(x))$ is the initial given velocity.
 The damping is from the resistance to the motion of the flow. It describes various physical situations such as porous media flow, drag or friction effects, and some dissipative mechanisms (see \cite{BD,BDC,H,HP} and references therein).
  The fact that ${\rm div}\,u = 0$, allows to write the term $(u.\nabla u):=u_1\partial_1 u+u_2\partial_2 u+u_3\partial_3u$ in the following form
$ {\rm div}\,(u\otimes u):=({\rm div}\,(u_1u),{\rm div}\,(u_2u),{\rm div}\,(u_3u)).$
 If the initial   velocity $u^0$ is quite regular, the divergence free condition determines the pressure $p$.\\
  Without loss of generality and in order to simplify the proofs of our results,  we consider the viscosity unitary ($\nu=1$).

  The global existence of weak solution of initial value problem of classical incompressible Navier-Stokes were proved by Leray and Hopf (see \cite{Hopf}-\cite{Leray}) long before. Uniqueness remains an open problem for the dimensions $d\geq3$.\\
  The polynomial damping $\alpha|u|^{\beta-1}u$ is studied in \cite{CJ} by Cai and Jiu, where they proved the global  existence of weak solution in
 $$L^\infty(\R^+,L^2(\R^3))\cap L^2(\R^+,\dot H^1(\R^3))\cap L^{\beta+1}(\R^+,L^{\beta+1}(\R^3)).$$
 The exponential damping $a (e^{b |u|^2}-1)u$ is studied in \cite{J1}  by J. Benameur, where he proved the global  existence of weak solution in
 $$L^\infty(\R^+,L^2(\R^3))\cap L^2(\R^+,\dot H^1(\R^3))\cap \mathcal E_b,$$
 where
 $\ds \mathcal E_b=\{f:\R^+\times\R^3\rightarrow\R\ :{\rm measurable}, \ (e^{b|f|^2}-1)|f|^2\in L^1(\R^+\times\R^3)\}.$\\

 The purpose of this paper is to study the well-posedness   and the asymptotic study   of the incompressible Navier-Stokes equations with exponential damping $ a (e^{b|u|^4}-1)u$.
We will show that the Cauchy problem $(S)$ has   a  global weak solutions for any $ a,b\in(0,\infty)$. We apply the Friedritch method to construct the approximate solutions and make more delicate estimates
to proceed to compactness arguments. In particular, we obtain new more a priory estimates:
$$\|u(t)\|_{L^2}^2+2\int_0^t\|\nabla u(s)\|_{L^2}^2ds+2 a\int_0^t\|(e^{b |u(s)|^4}-1)|u(s)|^2\|_{L^1}ds\leq \|u^0\|_{L^2}^2,$$
comparing with the
Navier-Stokes equations, to guarantee that the solution $u$ belongs to
$$L^\infty(\R^+,L^2(\R^3))\cap L^2(\R^+,\dot H^1(\R^3))\cap \mathcal F_b,$$
where  $\mathcal F_b=\{f:\R^+\times\R^3\to \R\ {\rm measurable},\  (e^{b |f|^4}-1)|f|^2\in L^1(\R^+\times\R^3)\}.$\\

  To prove the uniqueness we use an energy method and the approximates systems. The proof of the asymptotic study is based on a decomposition of the solution in high and low frequencies and the uniqueness of such  solution in a well chosen time $t_0$. 

In our case of exponential damping, we are trying to find more regularity of Leray
solution in $\cap_pL^p(\R^+,L^p(\R^3))$. In particular, we give a new energy estimate. Our main result is the following:

\begin{theorem}\label{th1}\pn
 Let $u^0\in L^2(\mathbb R^3)$ be a divergence free vector fields, then there is a unique global solution of the system  $(S)$:
$u\in C_b(\R^+,L^2(\mathbb R^3))\cap L^2(\R^+,\dot H^1(\mathbb
R^3))\cap\mathcal E_b$. Moreover, for all $t\geq0$
\begin{equation}\label{eqth1-1}
\|u(t)\|_{L^2}^2+2\int_0^t\|\nabla
u(s)\|_{L^2}^2ds+2a\int_0^t\|(e^{b
|u(s)|^4}-1)|u(s)|^2\|_{L^1}ds\leq \|u^0\|_{L^2}^2.
\end{equation}
Moreover, we have
\begin{equation}\label{eqth1-2}
\limsup_{t\to \infty}\|u(t)\|_{L^2}=0.
\end{equation}
\end{theorem}
\begin{rem}\pn
\begin{enumerate}
\item The new results in this theorem is the uniqueness of the global solution, the
continuity of the solution in the $L^2(\R^3)$  space and the asymptotic behavior at infinity.

\item Generally, for {\bf r}$\geq1$ the following
problem
\begin{equation}\label{$S_r$}
 \left\{ \begin{matrix}
     \partial_t u
 -\nu\Delta u+ u.\nabla u  +a (e^{b |u|^{\bf r}}-1)u =\;\;-\nabla p \hfill&{\hbox{\rm in }} \mathbb R^+\times \mathbb R^3\\
     {\rm div}\, u = 0 \hfill&{\hbox{\rm in }} \mathbb R^+\times \mathbb R^3\\
    u(0,x) =u^0(x) \;\;\hfill&{\hbox{\rm in }}\mathbb R^3,\\
    a,b>0\hfill&
\end{matrix}\right. \tag{$P_{\bf r}$}
\end{equation}
and by adapting the same proof of result of \cite{J1}, we show the
global existence of such a solution in $C_b(\R^+,L^2(\mathbb
R^3))\cap L^2(\R^+,\dot H^1(\mathbb R^3))\cap\mathcal E_b^{\bf
r}$, where
$$\mathcal E_b^{\bf r}=\{f:\R^+\times\R^3\to \R\;{\rm measurable};\;(e^{b |f|^{\bf r}}-1)|f|^2\in L^1(\R^+\times\R^3)\}.$$
Moreover, we get $$\|u(t)\|_{L^2}^2+2\int_0^t\|\nabla
u(s)\|_{L^2}^2ds+2a\int_0^t\|(e^{b |u(s)|^{\bf r}}-1)|u(s)|^2\|_{L^1}ds\leq \|u^0\|_{L^2}^2.$$ The asymptotic
result (\ref{eqth1-2}) is true for all {\bf r}$\geq\frac{7}3$, and the
index $\frac{7}3$ is a critical technical condition(See
(\ref{eqasym1})-(\ref{eqasym2})).
\end{enumerate}
 \end{rem}
\section{\bf Notations and Preliminary Results}
For a function  $f\colon\R^3\to\bar\R$    and $R>0$, the
Friedrich operator $J_R$  is defined by:
  $\ds J_R(D)f=\F^{-1}(\chi_{B_R} \widehat{f}),$
where $B_R$ is the ball of center $0$ and radius $R$. If $L^2_\sigma(\R^3)$ denotes the space of divergence-free vector fields in $L^2 (\R^3)$, the
 Leray projector $\mathbb P\colon (L^2(\R^3))^3\to (L^2_\sigma(\R^3))^3$ is  defined by:
$$\mathcal F(\mathbb P f)=\widehat{f}(\xi)-(\widehat{f}(\xi).\frac{\xi}{|\xi|})\frac{\xi}{|\xi|}=M(\xi)\widehat{f}(\xi),$$
where   $M(\xi)$ is the matrix $(\delta_{k,\ell}-\frac{\xi_k\xi_\ell}{|\xi|^2})_{1\leq k,\ell\leq 3} $.
If $u \in \mathcal S(\R^3)^3$,\\
$\ds \mathbb P(  u)_k(x) = \frac{1}{(2\pi)^{\frac 3 2}} \int_{\R^3}  \left( \delta_{kj}-\frac{\xi_k \xi_j}{ \vert \xi \vert^2}\right)
\widehat{  u}_j(\xi) \, e^{i \xi \cdot  x}\,   d\xi,$
where $\mathcal S(\R^n)$ is the  Schwartz space.
Define also the operator $A_R(D)$ on $L^2(\R^3)$ by:
 $$\ds A_R(D)u=\mathbb P J_R(D)u=\mathcal F^{-1}(M(\xi)\chi_{B_R}(\xi)\widehat{u}).$$
%


 To simplify the exposition of the main result, we first collect some preliminary results  and we give some new technical lemmas.

\begin{prop}(\cite{HBAF})\label{prop1}\pn
 Let $H$ be a Hilbert space.
\begin{enumerate}
\item The unit ball is weakly compact, that is: if $(x_n)$ is a bounded sequence   in $H$, then there is a subsequence $(x_{\varphi(n)})$ such that
$$(x_{\varphi(n)}|y)\to  (x|y),\;\forall y\in H.$$

\item If $x\in H$ and $(x_n)$   a bounded sequence   in $H$ such that
$\ds\lim_{n\to+\infty}(x_n|y)=  (x|y)$, for all $y\in H,$
then $\|x\|\leq\ds \liminf_{n\to \infty}\|x_n\|.$

\item If $x\in H$ and $(x_n)$ is a bounded sequence   in $H$ such that\\
$\ds\lim_{n\to+\infty}(x_n|y)=  (x|y)$, for all $y\in H$
and
 $\limsup_{n\to \infty}\|x_n\|\leq \|x\|,$
then $\ds \lim_{n\to \infty}\|x_n-x\|=0.$
\end{enumerate}
\end{prop}
We recall   the following product law in the homogeneous Sobolev spaces:

\begin{lem}(\cite{JYC})\label{lem1}\pn
Let $s_1,\ s_2$ be two real numbers and $d\in\N$.

\begin{enumerate}
\item If $s_1<\frac d 2$\; and\; $s_1+s_2>0$, there exists a constant  $C_1=C_1(d,s_1,s_2)$, such that: if $f,g\in \dot{H}^{s_1}(\mathbb{R}^d)\cap \dot{H}^{s_2}(\mathbb{R}^d)$, then $f.g \in \dot{H}^{s_1+s_2-\frac{d}{2}}(\mathbb{R}^d)$ and
$$\|fg\|_{\dot{H}^{s_1+s_2-\frac{d}{2}}}\leq C_1 (\|f\|_{\dot{H}^{s_1}}\|g\|_{\dot{H}^{s_2}}+\|f\|_{\dot{H}^{s_2}}\|g\|_{\dot{H}^{s_1}}).$$

\item If $s_1,s_2<\frac d 2$\; and\; $s_1+s_2>0$ there exists a constant $C_2=C_2(d,s_1,s_2)$ such that: if $f \in \dot{H}^{s_1}(\mathbb{R}^d)$\; and\; $g\in\dot{H}^{s_2}(\mathbb{R}^d)$, then  $f.g \in \dot{H}^{s_1+s_2-\frac{d}{2}}(\mathbb{R}^d)$ and
$$\|fg\|_{\dot{H}^{s_1+s_2-\frac{d}{2}}}\leq C_2 \|f\|_{\dot{H}^{s_1}}\|g\|_{\dot{H}^{s_2}}.$$
\end{enumerate}
 \end{lem}

\begin{lem}\label{lem2}\pn
    Let $ \beta>0$ and $d\in\N$. Then, for all $x,y\in\R^d$, we have

    \begin{equation}\label{eqn-lem2-1}
    \langle |x|^{\beta}x-|y|^{\beta}y ,x-y\rangle\geq \frac{1}{2}(|x|^{\beta}+|y|^{\beta})|x-y|^{2},
    \end{equation}
and, for ${\bf r}>0$, we have
\begin{equation}\label{eqn-lem2-2}
\langle (e^{b |x|^{\bf r}}-1)x-(e^{b |y|^{\bf
r}}-1)y,x-y\rangle\geq \frac{1}{2}\Big((e^{b |x|^{\bf
r}}-1)+(e^{b |y|^{\bf r}}-1)\Big)|x-y|^{2}.
\end{equation}
\end{lem}
{\bf Proof.} \pn
 Suppose that $|x|>|y|>0$. For $u>v>0$, we have

\begin{equation}\label{eqn-lem2-3}
2\langle ux-vy, x-y\rangle-(u+v)|x-y|^{2}=(u-v)(|x|^2-|y|^2)\geq0.
\end{equation}
It suffices to take $u=|x|^\beta$ and $v=|y|^\beta$, we get the inequality \eqref{eqn-lem2-1}.

\noindent
Suppose that $|x|>|y|>0$. In use of the inequality
\eqref{eqn-lem2-3} with $\ds u=(e^{b |x|^{\bf r}}-1)$ and
$\ds v=(e^{b |y|^{\bf r}}-1)$, we get

\beqq
&&2\langle (e^{b |x|^{\bf r}}-1)x-(e^{b |y|^{\bf r}}-1)y,x-y\rangle-  \Big((e^{b |x|^{\bf r}}-1)+(e^{b |y|^{\bf r}}-1)\Big)|x-y|^{2}\\
&&\hskip 3.5cm =(e^{b |x|^{\bf r}}- e^{b |y|^{\bf
r}})|x-y|^{2}\ge0. \eeqq This proves the inequality
\eqref{eqn-lem2-2}.

 The following result is a generalization of Proposition 3.1 in \cite{J1}.
\begin{prop}\label{prop2}  \pn
Let $\nu_1,\nu_2,\nu_3\in[0,\infty)$, $r_1,r_2,r_3\in(0,\infty)$ and $f^0\in L^2_\sigma(\R^3)$. \\
For $n\in\N$, let $F_n:\R^+\times\R^3\to \R^3$ be a measurable function in $C^1(\R^+,L^2(\R^3))$ such that $$A_n(D)F_n=F_n,\;F_n(0,x) =A_n(D)f^0(x)$$ and

\begin{enumerate}
\item [(E1)]
$\ds  \partial_t F_n+\sum_{k=1}^3\nu_k|D_k|^{2r_k} F_n+ A_n(D){\rm div}\,(F_n\otimes F_n)+ A_n(D)h(|F_n|)F_n =0.$

\item [(E2)]
\beqq
&&\ds \|F_n(t,.)\|_{L^2}^2+2\sum_{k=1}^3\nu_k\int_0^t\||D_k|^{r_k} F_n(s,.)\|_{L^2}^2ds\\
&&\hskip 2cm +2 a\int_0^t\|h(|F_n(s,.)|)|F_n(s,.)|^2\|_{L^1}ds \leq \|f^0\|_{L^2}^2.
\eeqq
\end{enumerate}
 where  $\ds h(z)= a(e^{b z^{\bf r}}-1),$\, with ${\bf r}\geq1$ and $a,b  >0$.
  Then: for every $\varepsilon>0$ there is $\delta=\delta(\varepsilon,a,b,\nu_1,\nu_2,\nu_3,r_1,r_2,r_3,\|f^0\|_{L^2})>0$
  such that: for all $t_1,t_2\in\R^+$, we have

\begin{equation}\label{eqn-1}
\Big(|t_2-t_1|<\delta\Longrightarrow \|F_n(t_2)-F_n(t_1)\|_{H^{-s_0}}<\varepsilon\Big),\;\forall n\in\N,
\end{equation}
with $\ds  s_0\ge \max(3,2r_1,2r_2,2r_3).$
\end{prop}

{\bf Proof.}\pn
 Integrate $(E1)$    on the interval $[t_1,t_2]\subset\R^+$ and take the inner product in $H^{-s_0}$, we get
\beqq \|F_n(t_2)-F_n(t_1)\|_{H^{-s_0}}&\leq&
\int_{t_1}^{t_2} \sum_{k=1}^3\nu_k \||D_k|^{2r_k}F_n(t)\|_{H^{-s_0}}dt\\
&&+
\int_{t_1}^{t_2}\|A_n(D){\rm div}\,(F_n\otimes F_n)(t)\|_{H^{-s_0}}dt\\
&&+\int_{t_1}^{t_2} \|A_n(D)h(|F_n|)F_n(t)\|_{H^{-s_0}}dt.
 \eeqq
Let

$$\ds I_{1,n}(t_1,t_2)=\int_{t_1}^{t_2} \sum_{k=1}^3\nu_k \||D_k|^{2r_k}F_n(t)\|_{H^{-s_0}}dt,$$

$$I_{2,n}(t_1,t_2)=\int_{t_1}^{t_2}\|A_n(D){\rm div}\,(F_n\otimes F_n)(t)\|_{H^{-s_0}}dt,$$
and
$$I_{3,n}(t_1,t_2)=\int_{t_1}^{t_2} \|A_n(D)h(|F_n|)F_n(t)\|_{H^{-s_0}}dt.$$
We have:
\begin{eqnarray}\label{eqn-2}
I_{1,n}(t_1,t_2)&\leq& \sum_{k=1}^3\nu_k\int_{t_1}^{t_2}\|F_n(t)\|_{H^{2r_k-s_0}}dt\\
&\buildu{2r_k-s_2\le0}\under\leq & \left(\sum_{k=1}^3\nu_k\right) \int_{t_1}^{t_2}\|F_n(t)\|_{L^2}dt\nonumber\\
&\leq&  \left(\sum_{k=1}^3\nu_k\right) \| f^0\|_{L^2}(t_2-t_1).\nonumber
\end{eqnarray}

 \beqq
\ds  I_{2,n}(t_1,t_2)&=& \int_{t_1}^{t_2}\|A_n(D){\rm
div}\,(F_n\otimes F_n)(t)\|_{H^{-s_0}}dt
 \leq \ds \int_{t_1}^{t_2}\|{\rm div}\,(F_n\otimes F_n)(s)\|_{H^{-s_2}}dt \\
&\leq&\ds \int_{t_1}^{t_2}\|(F_n\otimes F_n)(t)\|_{H^{-s_0+1}}dt \le\ds \int_{t_1}^{t_2}\|(F_n\otimes F_n)(t)\|_{H^{-2}}dt\\
 \eeqq

Recall that   if $f\colon\R^3\to \R^3$ is an integrable function, then for all $s>\frac 32$,

\beqq
\|f\|_{H^{-s}}^2&=& \int_{\R^3}(1+|\xi|^2)^{-s}|\widehat{f}(\xi)|^2d\xi
 \leq  \Big(\int_{\R^3}(1+|\xi|^2)^{-s}d\xi\Big)\|\widehat{f}\|_{ \infty }^2\\
& \leq& \Big(\int_{\R^3}(1+|\xi|^2)^{-s}d\xi\Big)\|f\|_{L^1 }^2.
\eeqq
We deduce that
\begin{equation}\label{eqn-3}
\|f\|_{H^{-s}}^2 \le \Big(\int_{\R^3}(1+|\xi|^2)^{-s}d\xi\Big)\|f\|_{L^1 }^2
\end{equation}
and there   exists   $C >0$ such that

\begin{eqnarray}\label{eqn-4}
 \ds  I_{2,n}(t_1,t_2)
 &\leq&  \ds  C \int_{t_1}^{t_2}\|(F_n\otimes F_n)(t)\|_{L^1}dt
 \leq \ds  C \int_{t_1}^{t_2}\|F_n(t)\|_{L^2}^2dt \\
 &\leq &\ds C (t_2-t_1)\|f^0\|_{L^2}^2.\nonumber
\end{eqnarray}

\noindent
To estimate the integral  $I_{3,n}(t_1,t_2)$, consider for $R>1$ the sub-level sets:
$$X_n(R,t)=\{x\in\R^3: |F_n(t,x)|\le R\}.$$
We remark that    $$(e^{b|F_n(t,x)|^{\bf r}}-1)|F_n(t,x)|\le(\frac{e^{bR^{\bf r}}-1}{R})  |F_n(t,x)|^2,\;\forall x\in X_n(R,t).$$

 Let $M(R)=\frac{e^{bR^{\bf r}}-1}{R}.$ From \eqref{eqn-3}, there exists $C_1,C_2,C_3>0$ such that

\beqq I_{3,n}(t_1,t_2)&=&\int_{t_1}^{t_2}
\|A_n(D)h(|F_n|)F_n(t)\|_{H^{-s_0}}dt
 \le C_1\int_{t_1}^{t_2} \| h(|F_n|)F_n(t)\|_{L^{1}}dt\\
&\le&C_1\int_{t_1}^{t_2}\int_{X_n(R,t)} h(|F_n|)|F_n|dx dt+C_1\int_{t_1}^{t_2}\int_{X_n(R,t)^c}h(|F_n|)|F_n|dx dt\\
&\le&C_2\int_{t_1}^{t_2}\int_{X_n(R,t)}  |F_n|^2dx dt+\frac{C_1}R\int_{t_1}^{t_2}\int_{X_n(R,t)^c}h(|F_n|)|F_n|^2dx dt\\
&\le&C_2M(R)\int_{t_1}^{t_2}\|F_n(t)\|_2^2dt+\frac{C_3}R\int_{t_1}^{t_2}\|h(|F_n|)|F_n|^2\|_{L^1} dt\\
&\le&C_2M(R) \|f^0\|_2^2(t_2-t_1)+\frac{C_3}R\|f^0\|_2^2.
\eeqq
Hence
\begin{equation}\label{eqn-5}
  I_{3,n}(t_1,t_2)\le C_2M(R) \|f^0\|_2^2(t_2-t_1)+\frac{C_3}R\|f^0\|_2^2.
\end{equation}

Now using the inequalities \eqref{eqn-2}, \eqref{eqn-4} and \eqref{eqn-5}, for $\varepsilon>0$, consider $R$ such that $\ds \frac{C_3}R\|f^0\|_2^2<\frac{\varepsilon}{4}$   and
$$0<\delta< \min\left(\frac{\varepsilon}{4 ((\sum_{k=1}^3\nu_k ) \| f^0\|_{L^2}+1)},\frac{\varepsilon}{4(C\|f^0\|_2^2+1)},
\frac{\varepsilon}{4(C_2M(R) \|f^0\|_2^2+1)}. \right).$$
For such  $\delta$, we get \eqref{eqn-1}.

\section{\bf Proof of the Main Theorem \ref{th1}}\ \\
The proof  is given in four steps:

\subsection{Existence of  Week Solution}\ \\
In this step, we   build approximate solutions of the system $(S)$ inspired by the method used in  \cite{J1,JYC}, hence we construct a global solution. For this, consider the approximate system with   parameter $n\in\N$:
$$(S_n)
  \begin{cases}
     \partial_t u
 -\Delta J_nu+ J_n(J_nu.\nabla J_nu)  + a J_n[(e^{b |J_nu|^4}-1)J_nu] =\;\;-\nabla p_n\hbox{ in } \mathbb R^+\times \mathbb R^3\\
 p_n=(-\Delta)^{-1}\Big({\rm div}\,J_n(J_nu.\nabla J_nu)  + a {\rm div}\,J_n[(e^{b |J_nu|^2}-1)J_nu]\Big)\\
     {\rm div}\, u = 0 \hbox{ in } \mathbb R^+\times \mathbb R^3\\
    u(0,x) =J_nu^0(x) \;\;\hbox{ in }\mathbb R^3.
  \end{cases}
$$
  $J_n$ is the Friedritch operator defined in the second section.

\begin{enumerate}
\item[$\bullet$] By Cauchy-Lipschitz Theorem, we obtain a unique solution $u_n\in C^1(\R^+,L^2_\sigma(\R^3))$ of $(S_{2,n})$. Moreover, $J_nu_n=u_n$ such that

\begin{equation}\label{eqn-6}
\|u_n(t)\|_{L^2}^2+2\int_0^t\|\nabla
u_n\|_{L^2}^2+2a\int_0^t\|(e^{b
|u_n|^4}-1)|u_n|^2\|_{L^1}\leq \|u^0\|_{L^2}^2.
\end{equation}

\item[$\bullet$]
The sequence $(u_n)_n$ is bounded in $L^\infty(\R^+,L^2(\R^3))$ and on $L^2(\R^+,\dot H^{1}(\R^3)$. Using proposition \ref{prop2} and the  interpolation method, we deduce that the sequence $(u_n)_n$ is equicontinuous on $H^{-1}(\R^3)$.

\item[$\bullet$] Let $(T_q)_q $ be a strictly  increasing sequence such that  $\ds\lim_{q\to+\infty} T_q=\infty$. Consider a sequence
of functions $(\theta_q)_{q }$ in $C_0^\infty(\R^3)$  such that

$$\left\{\begin{array}{l}
\theta_q(x)=1,\ {\rm for}\  |x|\le   q+\frac{5}{4}\\
\theta_q(x)=0,\ {\rm for}\   |x|\ge  q+2 \\
0\leq \theta_q\leq 1.
\end{array}\right.$$
Using  \eqref{eqn-6}, the equicontinuity of the sequence $(u_n)_n$   on $H^{-1}(\R^3)$
 and classical argument by combining Ascoli's theorem and the Cantor diagonal process, there exists a subsequence $(u_{\varphi(n)})_n$   and\\
$u\in L^\infty(\R^+,L^2(\R^3))\cap C(\R^+,H^{-3}(\R^3))$ such that: for all $q\in\N$,

\begin{equation}\label{eqn-7}
\lim_{n\to \infty}\|\theta_q(u_{\varphi(n)}-u)\|_{L^\infty([0,T_q],H^{-4})}=0.
\end{equation}
In particular, the sequence $(u_{\varphi(n)}(t))_n$ converges weakly in $L^2(\R^3)$ to $u(t)$ for all
$t\geq0$.

\item[$\bullet$] Combining the above inequalities, we obtain:

\begin{equation}\label{eqn-8}
\|u(t)\|_{L^2}^2\!+\!2\int_0^t\!\|\nabla
u(s)\|_{L^2}^2ds\!+\!2a\int_0^t\!\|(e^{b |u(s)|^4} -1)
|u(s)|^2\|_{L^1}ds\!\leq\! \|u^0\|_{L^2}^2.
\end{equation}
for all $t\geq0$.

\item[$\bullet$] $u$ is a solution of the system $(S)$.\end{enumerate}

\subsection{Continuity of the Solution in $L^2$}\pn
  In this section, we give a simple proof of the continuity of the solution $u$ of the system $(S)$ and we prove also that $u\in C(\R^+,L^2(\R^3))$. The construction of the solution is based on the Friedrich approximation method. We point out that we can use this method to show the same results as in  \cite{HP}.\\
$\bullet$ By inequality (\ref{eqn-8}), we get
$$\limsup_{t\to 0}\|u(t)\|_{L^2}\leq\|u^0\|_{L^2}.$$
Then, proposition \ref{prop1}-(3) implies that

$$\limsup_{t\to 0}\|u(t)-u^0\|_{L^2}=0,$$
which ensures the continuity of $u$ at $0$.\\
$\bullet$ Consider the  functions
$$v_{n,\varepsilon}(t,.)=u_{\varphi(n)}(t+\varepsilon,.),\;p_{n,\varepsilon}(t,.)=p_{\varphi(n)}(t+\varepsilon,.),$$
for $n\in\N$ and $\varepsilon>0$. We have:

$$\begin{array}{lcl}
\partial_tu_{\varphi(n)}-\Delta u_{\varphi(n)}+J_{\varphi(n)}(u_{\varphi(n)}.\nabla
u_{\varphi(n)})+a J_{\varphi(n)}(e^{b|u_{\varphi(n)}|^4}-1)u_{\varphi(n)}&=&-\nabla p_{\varphi(n)} \\
\partial_tv_{n,\varepsilon}-\Delta v_{n,\varepsilon}+J_{\varphi(n)}
(v_{n,\varepsilon}.\nabla v_{n,\varepsilon})+a J_{\varphi(n)}(e^{b|v_{n,\varepsilon}|^4}-1)v_{n,\varepsilon}&=&-\nabla p_{n,\varepsilon}\\
\end{array}$$
The function $w_{n,\varepsilon}=u_{\varphi(n)}-v_{n,\varepsilon}$ fulfills the following:

\beqq
&&\partial_tw_{n,\varepsilon}-\Delta w_{n,\varepsilon}  +a
J_{\varphi(n)}\Big((e^{b|u_{\varphi(n)}|^4}-1)u_{\varphi(n)}-(e^{b|v_{n,\varepsilon}|^4}-1)v_{n,\varepsilon}\Big)\\
&&\hskip  3cm = -\nabla (p_{\varphi(n)}-p_{n,\varepsilon})+J_{\varphi(n)}(w_{n,\varepsilon}.\nabla w_{n,\varepsilon})\\
&&\hskip  3cm-J_{\varphi(n)}(w_{n,\varepsilon}.\nabla u_{\varphi(n)})
 - J_{\varphi(n)}(u_{\varphi(n)}.\nabla w_{n,\varepsilon}).
 \eeqq
Taking the scalar product in $L^2(\R^3)$ with
$w_{n,\varepsilon}$ and using the properties ${\rm div}\ w_{n,\varepsilon}=0$ and
$\langle w_{n,\varepsilon}.\nabla w_{n,\varepsilon},w_{n,\varepsilon}\rangle=0$, we get

\begin{eqnarray}\label{eqn-9}
\frac{1}{2}\frac{d}{dt}\|w_{n,\varepsilon}\|_{L^2}^2+\|\nabla
w_{n,\varepsilon}\|_{L^2}^2& +&a \langle
J_{\varphi(n)}\Big((e^{b|u_{\varphi(n)}|^4}-1)u_{\varphi(n)}
-(e^{b|v_{n,\varepsilon}|^4}-1)v_{n,\varepsilon}\Big);w_{n,\varepsilon}\rangle_{L^2}
\nonumber\\
&=& -\langle J_{\varphi(n)}(w_{n,\varepsilon}.\nabla u_{\varphi(n)});w_{n,\varepsilon}\rangle _{L^2} .
\end{eqnarray}
Using inequality  \eqref{eqn-lem2-2}, we get

\beqq &&\langle
J_{\varphi(n)}\Big((e^{b|u_{\varphi(n)}|^2}-1)u_{\varphi(n)}-(e^{b|v_{n,\varepsilon}|^4}-1)v_{n,\varepsilon}\Big);
w_{n,\varepsilon}\rangle _{L^2}\\
 &&\hskip 5cm=\langle (e^{b|u_{\varphi(n)}|^4}-1)u_{\varphi(n)}-
 (e^{b|v_{n,\varepsilon}|^4}-1)v_{n,\varepsilon};J_{\varphi(n)}w_{n,\varepsilon}\rangle _{L^2}\\
 &&\hskip 5cm= \langle (e^{b|u_{\varphi(n)}|^4}-1)u_{\varphi(n)}-(e^{b|v_{n,\varepsilon}|^4}-1)v_{n,\varepsilon};w_{n,
\varepsilon}\rangle_{L^2}\\
 && \hskip 5cm\geq  \frac{1}{2}\int_{\R^3}\Big((e^{b|u_{\varphi(n)}|^4}-1)+(e^{b|v_{n,\varepsilon}|^4}-1)\Big)|w_{n,\varepsilon}|^2\\
&&\hskip 5cm  \geq   \frac{1}{2}\int_{\R^3}(e^{b|u_{\varphi(n)}|^4}-1)|w_{n,\varepsilon}|^2\\
 && \hskip 5cm\geq    \frac{b}{2}\int_{\R^3}|u_{\varphi(n)}|^4|w_{n,\varepsilon}|^2.
\eeqq

$$\begin{array}{lcl}
|\langle J_{\varphi(n)}(w_{n,\varepsilon}.\nabla u_{\varphi(n)});
w_{n,\varepsilon}\rangle _{L^2}|&\leq&\ds \int_{\R^3}|w_{n,\varepsilon}|.|u_{\varphi(n)}|.|\nabla w_{n,\varepsilon}|\\
&\leq&\ds \frac{1}{2}\int_{\R^3}|w_{n,\varepsilon}|^2|u_{\varphi(n)}|^2+\frac{1}{2}\|\nabla w_{n,\varepsilon}\|_{L^2}^2.
\end{array}$$
Again by using the elementary inequality $\ds xy\le
\frac{ab}8x^2+\frac{2}{ab}y^2$, for $x,y\ge 0$, we get

$$|\langle J_{\varphi(n)}(w_{n,\varepsilon}.\nabla u_{\varphi(n)});w_{n,\varepsilon}\rangle _{L^2}|
\leq\ds
\frac{ab}{8}\int_{\R^3}|u_{\varphi(n)}|^4|w_{n,\varepsilon}|^2
+\frac{2}{ab}\|w_{n,\varepsilon}\|_{L^2}^2+\frac{1}{2}\|\nabla
w_{n,\varepsilon}\|_{L^2}^2.$$
Combining the  identity  \eqref{eqn-lem2-2} and the inequality \eqref{eqn-9}, we get

$$\frac{1}{2}\frac{d}{dt}\|w_{n,\varepsilon}\|_{L^2}^2+\frac{1}{2}\|\nabla w_{n,\varepsilon}\|_{L^2}^2\leq
\frac{2}{ab}\|w_{n,\varepsilon}\|_{L^2}^2.$$
By Gronwall Lemma, we get
$$\|w_{n,\varepsilon}(t)\|_{L^2}^2\leq  \|w_{n,\varepsilon}(0)\|_{L^2}^2e^{\frac{4t}{ab} }.$$
But
$$\|u_{\varphi(n)}(t+\varepsilon)-u_{\varphi(n)}(t)\|_{L^2}^2\leq
\|u_{\varphi(n)}(\varepsilon)-u_{\varphi(n)}(0)\|_{L^2}^2e^{\frac{4t}{ab}
}.$$ For $t_0>0$ and $\varepsilon\in(0,t_0)$, we have

$$\|u_{\varphi(n)}(t_0+\varepsilon)-u_{\varphi(n)}(t_0)\|^{2}_{L^{2}}
\leq\|u_{\varphi(n)}(\varepsilon)-u_{\varphi(n)}(0)\|^{2}_{L^{2}}
\exp\Big(\frac{4t_0}{ab}\Big),$$

$$\|u_{\varphi(n)}(t_0-\varepsilon)-u_{\varphi(n)}(t_0)\|^{2}_{L^{2}}\leq
\|u_{\varphi(n)}(\varepsilon)-u_{\varphi(n)}(0)\|^{2}_{L^{2}}
\exp\Big(\frac{4t_0}{ab}\Big).$$
So
\beqq
\| u_{\varphi(n)}(\varepsilon)-u_{\varphi(n)}(0) \|_{L^2}^2 &=&
\| J_{\varphi(n)} u_{\varphi(n)}(\varepsilon)-J_{\varphi(n)}u_{\varphi(n)}(0) \|_{L^2}^2 \\
&=&\|\chi_{\varphi(n)}\left(\widehat{u_{\varphi(n)}}  -\widehat{u^0}\right)\|_{\varphi(n)}^2\\
&\le&\| u_{\varphi(n)}(\varepsilon)-u^0 \|_{L^2}^2\\
&\le& 2\|u^0 \|_{L^2}^2-2Re\langle
u_{\varphi(n)}(\varepsilon),u^0\rangle. \eeqq But $\ds
\lim_{n\to+\infty}\langle
u_{\varphi(n)}(\varepsilon),u^0\rangle=\langle u
(\varepsilon),u^0\rangle$. Hence
$$\liminf_{n\to \infty}\|u_{\varphi(n)}(\varepsilon)-u_{\varphi(n)}(0)\|^{2}_{L^{2}}\leq 2\|u^0\|^{2}_{L^{2}}-
2Re\langle u(\varepsilon);u^0\rangle_{L^2}.$$ Moreover, for all
$q,N\in\N$ \beqq
\|J_N\Big(\theta_q.(u_{\varphi(n)}(t_0\pm\varepsilon)-u_{\varphi(n)}(t_0))\Big)\|^{2}_{L^2}
 &\leq & \|\theta_q.(u_{\varphi(n)}(t_0\pm\varepsilon)-u_{\varphi(n)}(t_0))\|^{2}_{L^2}\\
 &\leq&
 \|u_{\varphi(n)}(t_0\pm\varepsilon)-u_{\varphi(n)}(t_0)\|^{2}_{L^2}.
 \eeqq

 Using  \eqref{eqn-7} we get, for $q$ big enough,
 $$\|J_N\Big(\theta_q.(u(t_0\pm\varepsilon)-u(t_0))\Big)\|^{2}_{L^2}
 \leq \liminf_{n\to \infty}\|u_{\varphi(n)}(t_0\pm\varepsilon)-u_{\varphi(n)}(t_0)\|^{2}_{L^2}.$$
 Then
$$\|J_N\Big(\theta_q.(u(t_0\pm\varepsilon)-u(t_0))\Big)\|^{2}_{L^2}
 \leq 2\Big(\|u^0\|^{2}_{L^{2}}-Re\langle u(\varepsilon);u^0\rangle_{L^2}\Big)\exp(\frac{4t_0}{ab}).$$
By applying the monotone convergence theorem in the order $N\to \infty$ and $q\to \infty$, we get
$$\|u(t_0\pm\varepsilon,.)-u(t_0,.))\|^{2}_{L^2}
 \leq 2\Big(\|u^0\|^{2}_{L^{2}}-Re\langle u(\varepsilon);u^0\rangle_{L^2}\Big)\exp(\frac{4t_0}{ab}).$$
Using the continuity at 0 and make $\varepsilon\to 0$, we get the continuity at $t_0$.

\subsection{Uniqueness of the Solution}\ \\
Let $u,v$ be two solutions of $(S)$ in the space

$$C_b(\R^+,L^2(\R^3))\cap L^2(\R^+,\dot H^1(\R^3))\cap \mathcal F_\beta.$$
The function $w=u-v$ satisfies the following:

$$\partial_tw-\Delta w+a \Big((e^{b|u|^4}-1)u-(e^{b|v|^4}-1)v\Big)= -\nabla (p-\tilde p)+w.\nabla w-w.\nabla u-
u.\nabla w  .$$
Taking the scalar product in $L^2$ with $w$, we get

$$\frac{1}{2}\frac{d}{dt}\|w\|_{L^2}^2+\|\nabla w\|_{L^2}^2+a \langle \Big((e^{b|u|^4}-1)u-(e^{b|v|^2}-1)v\Big);w\rangle _{L^2}=-\langle w.\nabla u;w\rangle _{L^2}  .$$
The idea is to lower the term $\langle \Big((e^{b|u|^4}-1)u-(e^{b|v|^2}-1)v\Big);w\rangle _{L^2}$ with the help of the Lemma \ref{lem2} and then divide the term find into two equal pieces, one to absorb the nonlinear term and the other is used in the last inequality.\\

By using inequality \eqref{eqn-lem2-2}, we get

\beqq
\langle \Big((e^{b|u|^4}-1)u-(e^{b|v|^4}-1)v\Big);w\rangle _{L^2}&\geq & \frac{1}{2}\int_{\R^3}  \Big((e^{b|u|^4}-1)+(e^{b|v|^4}-1)\Big)|w|^2\\
&\geq& \frac{b}{2}\int_{\R^3}|u|^4|w|^2.
\eeqq
Moreover, we have
\beqq
|\langle w.\nabla u;w\rangle _{L^2}|&=&|\langle {\rm div}\,(w\otimes u);w\rangle _{L^2}|=|\langle w\otimes u;\nabla w\rangle _{L^2}|\\
&\leq&\ds \int_{\R^3}|w|.|u|.|\nabla w|
\leq\ds \frac{1}{2}\int_{\R^3}|w|^2|u|^2+\frac{1}{2}\|\nabla w\|_{L^2}^2\\
&\leq&\ds  \frac{ab}{8}\int_{\R^3}|u|^4|w|^2+\frac{1}{2ab}\|w\|_{L^2}^2+\frac{1}{2}\|\nabla w\|_{L^2}^2\\
\eeqq
Combining the above inequalities, we get
$$\frac{1}{2}\frac{d}{dt}\|w\|_{L^2}^2+\frac{1}{2}\|\nabla w\|_{L^2}^2+\frac{a}{4}\int_{\R^3}  \Big((e^{b|u|^2}-1)+(e^{b|v|^2}-1)\Big)|w|^2\leq \frac{1}{2ab}\|w\|_{L^2}^2$$
and, Gronwall Lemma gives
$$\|w\|_{L^2}^2+\int_0^t\|\nabla w\|_{L^2}^2+\frac{a}{2}\int_0^t\int_{\R^3}  \Big((e^{b|u|^4}-1)+(e^{b|v|^4}-1)\Big)|w|^2\leq \|w^0\|_{L^2}^2e^{\frac{t}{ab}}.$$
As $w^0=0$, then $w=0$ and $u=v$. Which implies the uniqueness.

\subsection{Asymptotic Study of the Global Solution}\pn

In this subsection we  prove the asymptotic behavior (\ref{eqth1-2}). For this we prove some preliminaries lemmas:

\begin{lem}\label{lem1}\pn
 If $u$ is a global solution of \eqref{$S_2$}, then  $(e^{b|u|^4}-1)u\in L^1(\R^+\times\R^3)$.
\end{lem}
\noindent{\bf Proof}\pn
If $X_1=\{(t,x):\ |u(t,x)|\leq 1\}$ and $X_2=\{(t,x):\ |u(t,x)|> 1\},$
we  have
$$\int_0^\infty\|(e^{b|u(s,.)|^4}-1)u)(s,.)\|_{L^1}ds=K_1+K_2$$
where

$$ K_1=\ds\int_{X_1}(e^{b|u(s,x)|^4}-1)|u(s,x)|dxds \quad {\rm and}\quad
K_2=\ds\int_{X_2}(e^{b|u(s,x)|^4}-1)|u(s,x)|dxds.$$

 \beqq
K_1&=&\ds\int_{X_1}(e^{b|u(s,x)|^4}-1)|u(s,x)|dxds \\
&=&\ds\int_{X_1}b|u(s,x)|^{\frac{5}3}\frac{(e^{b|u(s,x)|^4}-1)}{b|u(s,x)|^{4}}|u(s,x)|^{\frac{10}3}dxds \\
&\leq&\ds b(e^{b}-1)\int_{0}^\infty\int_{\R^3}|u(s,x)|^{\frac{10}3}dxds
 \leq \ds be^{b}\int_{0}^\infty\|u(s,.)\|_{L^{\frac{10}3}}^{\frac{10}3}ds .
\eeqq

By using the Sobolev injection
$\dot H^{\frac{3}5}(\R^3)\hookrightarrow L^{\frac{10}3}(\R^3)$, we get
\begin{equation}\label{eqasym1}
K_1\leq C \int_{0}^\infty\|u(s,.)\|_{\dot H^{\frac{3}5}}^{\frac{10}3}ds.
\end{equation}

By interpolation inequality
$\ds \|u(s)\|_{\dot H^{\frac{3}5}}\leq \|u(s)\|_{\dot H^0}^{\frac{2}5}\|u(s)\|_{\dot H^1}^{\frac{3}5}$,
we obtain
\begin{equation}\label{eqasym2}
K_1\leq C \int_{0}^\infty\|u(s,.)\|_{L^2}^{\frac{4}3}\|\nabla u(s)\|_{L^2}^2\leq C
\|u^0\|_{L^2}^{\frac{4}3}\int_{0}^\infty\|\nabla u(s)\|_{L^2}^2.\end{equation}
For the therm $K_2$, we have

$$K_2=\ds\int_{X_2}(e^{b|u(s,x)|^4}-1)|u(s,x)|dxds\le \ds\int_{0}^\infty\int_{\R^3}(e^{b|u(s,x)|^4}-1)|u(s,x)|^2dxds.
$$

Hence
$$\|(e^{b|u|^4}-1)u\|_{L^1(\R^+\times\R^3)}\leq C
\|u^0\|_{L^2}^{\frac{4}3}\int_{0}^\infty\|\nabla
u(s,.)\|_{L^2}^2ds +\int_{0}^\infty\int_{\R^3}(e^{b|u(s,x)|^4}-1)|u(s,x)|^2dxds.$$
Therefore $(e^{b|u|^4}-1)u\in L^1(\R^+\times\R^3)$.

\begin{lem}\pn
 If $u$ is a global solution of \eqref{$S_2$}, then  $\lim_{t\rightarrow\infty}\|u(t)\|_{H^{-2}}=0$.
\end{lem}
\noindent{\bf Proof}\pn
Let $\varepsilon>0$. By the energy inequality
(\ref{eqth1-1}) and Lemma \ref{lem1}, there exists   $t_0\geq0$ such that

\begin{equation}\label{asym.eq1}
\|\nabla u\|_{L^2([t_0,\infty)\times\R^3)}<\frac{\varepsilon}{4},
\end{equation}

\begin{equation}\label{asym.eq2}
\|(e^{b|u|^4}-1)u\|_{L^1([t_0,\infty)\times\R^3)}<\frac{\varepsilon}{4}.
\end{equation}
Now, consider the following system
\begin{equation}\label{$4.6$}
 \left\{ \begin{matrix}
     \partial_t v
 -\nu\Delta v+ v.\nabla v  +a (e^{b |v|^4}-1)v =\;\;-\nabla q \hfill&\hbox{ in } \mathbb R^+\times \mathbb R^3\\
     {\rm div}\, v= 0 \hfill&\hbox{ in } \mathbb R^+\times \mathbb R^3\\
    v(0,x) =u(t_0,x) \;\;\hfill&\hbox{ in }\mathbb R^3\hfill&.
\end{matrix}\right. \tag{$S'$}
\end{equation}
By the existence and uniqueness part, the system ($S'$) has a
unique global solution $v$ such that $v(t_0)=u(t_0,x)$ and $q(t)=p(t_0+t).$
We   recall  the following energy estimate  for this system:
$$\|v(t)\|_{L^2}^2+2\int_0^t\|\nabla
v(s)\|_{L^2}^2ds +2a\int_0^t\|(e^{b|v(s)|^2}-1)|v(s)|^2\|_{L^1}\leq\|u(t_0)\|_{L^2}^2\leq
\|u^0\|_{L^2}^2.$$
By the Duhamel formula,   $\ds v(t,x)=e^{t\Delta}v^0(x)+f(t,x)+g(t,x)$, where
$$f(t,x)=-\int_0^te^{(t-s)\Delta}\mathbb P{\rm div\,}(v\otimes v)(s,x)ds$$
and
$$g(t,x)=-\int_0^te^{(t-s)\Delta}\mathbb P{\rm div\,}(e^{b|v(s,x)|^2}-1)v(s,x)ds.$$
By the  Dominated Convergence Theorem, we have:
 $\ds \lim_{t\rightarrow\infty}\|e^{t\Delta}v^0\|_{L^2}=0$ and hence $\ds  \lim_{t\rightarrow\infty}\|e^{t\Delta}v^0\|_{H^{-2}}=0.$
 Moreover, we have
\beqq
\|f(t)\|_{H^{-2}}^2&\leq&\ds\|f(t)\|_{H^{-\frac{1}2}}^2 \leq\ds\|f(t)\|_{\dot H^{-\frac{1}2}}^2\\
&\leq&\ds\int_{\R^3}|\xi|^{-1}\left(\int_0^te^{-(t-s)|\xi|^2}|\mathcal F{\rm div}(v\otimes v)(s,\xi)|ds\right)^2d\xi\\
&\leq&\ds\int_{\R^3}|\xi|\left(\int_0^te^{-(t-s)|\xi|^2}|\mathcal F(v\otimes v)(s,\xi)|ds\right)^2d\xi.
\eeqq
 As
\beqq
\ds\left(\int_0^te^{-(t-s)|\xi|^2}|\mathcal F(v\otimes v)(s,\xi)|ds\right)^2
&\leq&\ds\left(\int_0^te^{-2(t-s)|\xi|^2}ds\right) \int_0^t|\mathcal F(v\otimes v)(s,\xi)|^2ds \\
&\leq&\ds|\xi|^{-2} \int_0^t|\mathcal F(v\otimes v)(s,\xi)|^2ds,
\eeqq
we obtain
\beqq
\|f(t)\|_{H^{-2}}^2dt&\leq&\ds\int_{\R^3}|\xi|^{-1}\int_0^t|\mathcal F(v\otimes v)(s,\xi)|^2dsd\xi \\
&\leq&\ds\int_0^t(\int_{\R^3}|\xi|^{-1}|  (v\otimes v)(s,\xi)|^2d\xi)ds=\ds\int_0^t\|v\otimes v)(s)\|_{\dot H^{-\frac{1}2}}^2ds.
\eeqq

\noindent
Using the  product law in homogeneous Sobolev spaces, with
$s_1=0,\,s_2=1$, we get
\beqq
\|f(t)\|_{H^{-2}}^2dt&\leq&\ds C\int_0^t\|v(s)\|_{L^2}^2\|\nabla
v(s)\|_{L^2}^2ds.
\eeqq
 Using inequalities \eqref{asym.eq1}  and \eqref{asym.eq2}, we get
\beqq
\|f(t)\|_{H^{-2}}^2dt&\leq&\ds C\|u^0\|_{L^2}^2\int_0^t\|\nabla
u(t_0+s)\|_{L^2}^2ds\\
&\leq&\ds C\|u^0\|_{L^2}^2\int_0^\infty\|\nabla u(t_0+s)\|_{L^2}^2ds\\
&\leq&\ds C\|u^0\|_{L^2}^2\int_{t_0}^\infty\|\nabla u(s)\|_{L^2}^2ds\\
&\leq&\ds C\|u^0\|_{L^2}^2\frac{\varepsilon^2}{9(C\|u^0\|_{L^2}^2+1)},
\eeqq
which implies
$$\|f(t)\|_{H^{-2}}<\frac \varepsilon3,\;\forall t\geq 0.$$

For an estimation of $\|g(t)\|_{H^{-2}}$ and by using (\ref{eqn-3}) with $s=2$, we get
\beqq
\|g(t)\|_{H^{-2}}^2dt&\leq&\ds
\int_{\R^3}(1+|\xi|^2)^{-2}\left(\int_0^te^{-(t-s)|\xi|^2}|\mathcal
F((e^{b|v |^4}-1)v )(s,\xi)|ds\right)^2d\xi\\
&\leq&\ds
C\left(\int_0^t\|(e^{b|v(s,.)|^4}-1)v(s)\|_{L^1(\R^3)}ds\right)^2\\
&\leq&\ds C\|(e^{b|v|^4}-1)v\|_{L^1(\R^+\times\R^3)}^2,
\eeqq
where $\ds C=\int_{\R^3}(1+|\xi|^2)^{-2}d\xi.$
 Also by using inequality (\ref{asym.eq2}), we get
\beqq
\|g(t)\|_{H^{-2}}^2dt&\leq&C\|(e^{b|u(t_0+.)|^4}-1)u(t_0+.)\|_{L^1(\R^+\times\R^3)}^2\\
&\leq&C\|(e^{b|u|^4}-1)u\|_{L^1([t_0,\infty)\times\R^3)}^2
 \leq C\frac{\varepsilon^2}{9C}.
\eeqq

which implies that
$\ds \|g(t)\|_{H^{-2}}<\frac \varepsilon 3,\;\forall t\geq 0.$

Combining the above inequalities, we obtain

$$\lim_{t\rightarrow\infty}\|u(t)\|_{H^{-2}}=0.$$

\begin{lem}\pn
 If $u$ is a global solution of \eqref{$S_2$}, then  $\lim_{t\rightarrow\infty}\|u(t)\|_{L^2}=0.$
\end{lem}
\noindent{\bf Proof}\pn
 We have $u=w_1+w_2$, where
$$w_1 = {\bf 1}_{|D|<1}u=\mathcal F^{-1}\big({\bf
1}_{|\xi|<1}\widehat{u}\big),\qquad
w_2 = {\bf 1}_{|D|\geq1}u=\mathcal F^{-1}\big({\bf
1}_{|\xi|\geq1}\widehat{u}\big).$$
By the second step, we get
$$\|w_1(t)\|_{L^2}=c_0\|w_1(t)\|_{H^0}\leq 2c_0\||w_1(t)\|_{H^{-2}}\leq 2\||u(t)\|_{H^{-2}},$$
which implies
$$\lim_{t\rightarrow\infty}\|w_1(t)\|_{L^2}=0.$$
Let $\varepsilon>0$. There is a time $t_1>0$ such that
$$\|w_1(t)\|_{L^2}<\frac \varepsilon2,\;\forall t\geq t_1.$$
We have
$$\int_{t_1}^\infty\|w_2(t)\|_{L^2}^2dt\leq \int_{t_1}^\infty\|\nabla w_2(t)\|_{L^2}^2dt\leq \int_{t_1}^\infty\|\nabla u(t)\|_{L^2}^2dt<\infty.$$
As $t\rightarrow\|w_2(t)\|_{L^2}$ is continuous, then there is a
time $t_2\geq t_1$ such that $$\|w_2(t_2)\|_{L^2}<\frac \varepsilon2.$$
Particularly
$$\|u(t_2)\|_{L^2}^2=\|w_1(t_2)\|_{L^2}^2+\|w_2(t_2)\|_{L^2}^2<\frac {\varepsilon^2}2.$$
By using the following energy estimate
$$
\|u(t)\|_{L^2}^2+2\int_{t_2}^t\|\nabla
u(s)\|_{L^2}^2ds+2a\int_{t_2}^t\|(e^{b|u(s)|^2}-1)|u(s)|^2\|_{L^1}ds\leq\|u(t_2)\|_{L^2}^2,\,\forall t\geq t_2,$$ we get
$$\|u(t)\|_{L^2}<\varepsilon,\;\forall t\geq t_2,$$
and the proof is completed.

\section*{Acknowledgement}
This Project was funded by the National Plan for Science, Technology and Innovation (Maarifah), King Abdulaziz City for Science and Technology, Kingdom of Saudi Arabia, Award Number (14-MAT730-02).

\end{document}